\documentstyle[twoside,12pt, oldlfont]{article}

\begin{document}
\title {Sharper changes in topologies }
\author{Greg Hjorth}
\date{\today }                       

\input{amssym.def}
\input{amssym.tex}
\def\Ubf#1{{\baselineskip=0pt\vtop{\hbox{$#1$}\hbox{$\sim$}}}{}}
\maketitle
\begin{abstract} Let $G$ be a Polish group, $\tau$ a Polish topology on a space $X$, 
$G$ acting continuously on $(X,\tau)$, with $B\subset X$ $G$-invariant and in the Borel 
algebra generated by $\tau$. Then there is a larger Polish topology $\tau^*\supset \tau$ 
on $X$ so that $B$ is open with respect to $\tau^*$, $G$ still acts 
continuously on $(X,\tau^*)$, {\it and} $\tau^*$ has a basis 
consisting of sets that are of the same Borel rank as $B$ relative to $\tau$. 
\end{abstract}
{\bf \S 0. The preface}\\

0.1. Theorem(classical): Let $(X, \tau)$ be a Polish space. Let $B\subset X$ be Borel 
with respect to $\tau$. Then there is a richer topology $\tau^*\supset\tau$ such that

(i) every $\tau^*$-open set is Borel in $(X,\tau)$; 

(ii) $(X,\tau^*)$ is still a Polish space; 

(iii) $B$ is open with respect to $\tau^*$. (See \cite{kechris}.)\\

Thus we have that from the point of view of properties that can be stated solely with 
reference to Borel structure, every theorem provable for open sets 
holds as well for arbitrary Borel sets. For instance, as shown in \cite{kechris}, 
we obtain a fast proof that every uncountable Borel set has size $2^{\aleph_0}$. 

It turns out that this classical theorem can also be proved in the dynamical context.\\

0.2. Theorem(Sami). Let $S_{\infty}$ be the permutation group on a countably infinite 
set, viewed as a topological group with the topology of pointwise convergence. 
Suppose $S_{\infty}$ acts continuously on a Polish space $(X,\tau)$ and that 
$B\subset X$ is Borel with respect to $\tau$ and $S_{\infty}$-invariant. 
Then there is a richer topology $\tau^*\supset\tau$ such that

(i) every $\tau^*$-open set is Borel in $(X,\tau)$;

(ii) $(X,\tau^*)$ is still a Polish space;

(iii) $B$ is open with respect to $\tau^*$; 

(iv) $S_{\infty}$ acts continuously on  $(X,\tau^*)$.(See \cite{sami}.)\\

Then later:\\

0.3. Theorem(Becker-Kechris). Let $G$ be a Polish group. 
Suppose $G$ acts continuously on a Polish space $(X,\tau)$ and that
$B\subset X$ is Borel with respect to $\tau$ and $G$-invariant.
Then there is a richer topology $\tau^*\supset\tau$ such that

(i) every $\tau^*$-open set is Borel in $(X,\tau)$;

(ii) $(X,\tau^*)$ is still a Polish space;

(iii) $B$ is open with respect to $\tau^*$;

(iv) $G$ acts continuously on  $(X,\tau^*)$.(See \cite{beckerkechris}.)\\

Thus it may seem that a happy story has come to a pleasing conclusion. 

Yet there remained a gap in our understanding of these theorems regarding changes in 
Polish topologies. Whereas with 0.1 and 0.2 it was shown that the change could be effected 
with the minimum possible disturbance to $\tau$, with 0.3 we had only a crude upper 
bound for the Borel complexity -- with respect to $\tau$ -- of the open sets in $\tau^*$. 
Thus the authors of \cite{beckerkechris} were led to ask whether we can in general 
choose $\tau^*$ so that its basic open sets have approximately the same Borel 
complexity as $B$. 

Here we show this. In particular:\\

0.4. Theorem. Let $G$ be a Polish group.
Suppose $G$ acts continuously on a Polish space $(X,\tau)$ and that
$B\subset X$ is $F_{\sigma}$ with respect to $\tau$ and $G$-invariant.
Then there is a richer topology $\tau^*\supset\tau$ such that

(i) every $\tau^*$-open set is $F_{\sigma}$ in $(X,\tau)$;

(ii) $(X,\tau^*)$ is still a Polish space;

(iii) $B$ is open with respect to $\tau^*$;

(iv) $G$ acts continuously on  $(X,\tau^*)$.\\

And more generally:\\ 

0.5. Theorem. Let $G$ be a Polish group and let $\alpha$ be a countable 
ordinal.
Suppose $G$ acts continuously on a Polish space $(X,\tau)$ and that
$B$ is $\Ubf{\Sigma}^0_{\alpha}(X,\tau)$ and $G$-invariant.
Then there is a richer topology $\tau^*\supset\tau$ such that

(i) every $\tau^*$-open set is $\Ubf{\Sigma}^0_{\alpha}(X,\tau)$;

(ii) $(X,\tau^*)$ is still a Polish space;

(iii) $B$ is open with respect to $\tau^*$;

(iv) $G$ acts continuously on  $(X,\tau^*)$.\\ 

This answers question 5.1.9 of \cite{beckerkechris}. 

An unexpected advantage of the new proof is its brevity. 

\newpage
{\bf \S1. The landscape}\\

1.1. Definition. A topological group is said to {\it Polish} if it is Polish as a 
topological space -- which is to say that it is separable and it allows a complete 
compatible metric. If $G$ is a Polish group acting continuously on a Polish space $X$, 
then I will say that $X$ is a {\it Polish} $G$-{\it space}.\\

1.2. Definition. Let $G$ be a Polish group and $X$ a Polish $G$-space. 
For $U\subset G$ and $B\subset X$ we define the {\it Vaught transforms} by 
the specification that $B^{\Delta U}$ is the set of $x\in X$ such that 
$\{g\in G: g\cdot x\in B\}$ is non-meager and that $B^{* U}$ is 
the set of $x\in X$ such that      
$\{g\in G: g\cdot x\in B\}$ is comeager.\\

1.3. Lemma(Vaught). Let $(A_i)$ be a sequence of Borel sets, and let ${\cal B}$ be a 
basis of $G$. Fix $U\subset G$ open, $U\neq\emptyset$. Then 
\[(\bigcup A_i)^{\Delta U}=\bigcup\{A_i^{*V}|V \subset U, V\neq\emptyset, V\in{\cal B}\}.\]

(See \cite{vaught}.)\\


As a word to notation, if $X$ is a set and $\tau$ is a topology on $X$, I will try to use 
$X$ to denote the topological space -- that is, $X$ equipped with the topology -- as long as 
the intention is clear. When more than one topology is being considered on $X$, we will need 
to be specific, and instead use $(X,\tau)$, or $(X,\tau^*)$, and so on.\\

1.4. Definition. If $G$ is a group and $d$ is a metric on $G$, then $d$ is said to be {\it 
right invariant} if for all $g, h_1, h_2 \in G$ 
\[d(h_1\cdot g,h_2\cdot g)=d(h_1,h_2).\]

One can similarly define the notion of left invariant metric. Note that these are companion 
notions, since if $d$ is left invariant, then $d^*$ defined by $d^*(g,h)=d(g^{-1},h^{-1})$ will 
be right invariant.\\

1.5. Theorem(Birkhoff-Kakutani). Any Polish group has a compatible right invariant metric. 
(See \cite{hewittross}, 8.3.)\\

The metric provided by 1.5 will be compatible with the topology, 
but not necessarily complete.\\ 

1.6. Corollary. Any Polish group has a compatible right invariant metric bounded by 1 -- that 
is $\forall g,h\in G$ $d(g,h)\leq 1$. 

Proof. Given $d$ as in 1.5, let $d^*(g,h)=d(g, h)/[1+d(g,h)]$.$\Box$\\

1.7. Definition. Let $Y$ be a topological space and let $d$ be a compatible metric bounded by 1. Then 
let ${\cal L}(Y,d)$ be the set of (necessarily continuous) functions 
\[f:Y\rightarrow [0,1]\] such that for all $y_1, y_2\in Y$
\[|f(y_1)-f(y_2)|\leq d(y_1,y_2).\]

1.8. Lemma. ${\cal L}(Y,d)$ is compact in the topology of pointwise 
convergence. If $Y$ is separable then it is also metrizable. 

Proof. The first statement is Tychonov's theorem, while the second follows since for 
$Q=(a_i)_{i\in{\Bbb N}}\subset 
Y$ dense and countable we can identify ${\cal L}(Y,d)$ with a closed subset 
of $[0,1]^Q$; this last space has a metric given by 
$d(\vec x,\vec y)=\Sigma\{2^{-i}\cdot |x(a_i)-y(a_i)|:
i \in{\Bbb N}\}$.$\Box$\\

1.9. Definition. Let $G$ be a Polish group and let $d$ be a right invariant compatible metric on 
$G$ that is bounded by 1. Then for $f\in{\cal L}(G,d)$ let $g\cdot f$ be defined by $(g\cdot f)(h)=
f(hg)$ for all $h\in G$.\\

1.10. Lemma. This defines an action on ${\cal L}(G,d)$ under which it becomes a compact 
Polish $G$-space. 

Proof. By right invariance of $d$, along with the observation that 
\[(g_1g_2\cdot f)(h)=f(hg_1g_2)=g_2\cdot f(hg_1)=(g_1\cdot(g_2\cdot f))(h),\]
this defines an action of $G$ on ${\cal L}(G,d)$. 
Continuity follows by the definition of the space and its topology,
since if we fix $h\in G$, and let $U$ be an open neighborhood of $h$ such that for all 
$h'\in U(d(h, h')<\epsilon)$,  
then for $g\in G$ small enough to ensure that $hg\in U$, and all  
$f\in {\cal L}(G,d)$ 
\[|g\cdot f(h)-f(h)|=_{df}|f(hg)-f(h)|\leq d(hg,h)<\epsilon.\] 
Compactness is 1.8.$\Box$\\

For the sake of being thorough:\\

1.11. Definition. Let $(X,\tau)$ be a Polish space (with $\tau$ the Polish topology). Then 
the $\Ubf{\Sigma}^0_1(X,\tau)$ sets are the open sets in $(X,\tau)$; for $\alpha$ a countable 
ordinal, a set is 
$\Ubf{\Pi}^0_{\alpha}(X,\tau)$ if its complement is $\Ubf{\Sigma}^0_{\alpha}(X,\tau)$;  
and given that we have defined $\Ubf{\Sigma}^0_{\alpha}(X,\tau)$ all $\alpha<\delta$, 
we say that a set is $\Ubf{\Sigma}^0_{\delta}(X,\tau)$ if it has the form $\bigcup \{A_i:i\in{\Bbb N}\}$ 
where each $A_i$ is $\Ubf{\Pi}^0_{\alpha(i)}(X,\tau)$ for some $\alpha(i)<\delta$.\\

Thus the $\Ubf{\Pi}^0_2(X,\tau)$ sets are precisely the $G_{\delta}$ sets in $(X,\tau)$, 
and the $\Ubf{\Sigma}^0_2(X,\tau)$ are the $F_{\sigma}$. Note by 1.3 and a transfinite induction, 
if $(X,\tau)$ is a Polish $G$-space, $B\in \Ubf{\Sigma}^0_{\alpha}(X,\tau)$ 
(or $\Ubf{\Pi}^0_{\alpha}(X,\tau)$), $U\subset G$ open, then $B^{\Delta U}$ (respectively 
$B^{*U}$) is  
$\Ubf{\Sigma}^0_{\alpha}(X,\tau)$ (or $\Ubf{\Pi}^0_{\alpha}(X,\tau)$ respectively). 
The base step of this induction begins with the observation that 
if $B$ is open, then $B^{\Delta U}=\{x\in X:\exists g\in U(g\cdot x \in B)\}$, and so is open.\\ 

1.12. Lemma(classical): Let $(X,\tau)$ be a Polish space and let $X_0$ be $\Ubf{\Pi}^0_2(X,\tau)$. 
Then $X_0$ is Polish in the relative topology. (See \cite{kechris}.)\\

1.13. Lemma(classical): Let $(X_i)_{i\in{\Bbb N}}$ be a sequence of Polish spaces; then 
$\Pi\{X_i:i\in{\Bbb N}\}$ in the product topology is Polish. (See \cite{kechris}, or use the 
proof of the second half of 1.8.)\\

1.14. Lemma(classical): Let $(\tau_i)_{i\in{\Bbb N}}$ be an increasing sequence of 
Polish topologies on $X$. Then the topology generated by the union $\bigcup \{\tau_i:i\in{\Bbb N}\}$ 
is Polish. (See \cite{kechris} or \cite{sami}.)\\

1.15. Definition. Let $G$ be a Polish group, $X$ and $Y$ Polish $G$-spaces. $\pi:X\rightarrow Y$ is 
a $G$-{\it mapping} if for all $g\in G$, $x\in X$
\[\pi(g\cdot x)=g\cdot \pi(x);\]
it is a {\it continuous} (or {\it Borel}) $G$-mapping if it is also continuous (respectively, 
the pullback of open 
sets are Borel); 
it is a $G$-{\it embedding} if it is one to one.\\

1.16. Definition. Let $d$ be a right invariant metric on $G$ and ${\cal O}\subset X$ open. 
Then $\varphi^{\cal O}:X\rightarrow {\cal L}(G, d)$ is the map $x\mapsto \varphi^{\cal O}_x$, 
where for $h\in G$
\[\varphi^{\cal O}_x(h)=d(h,\{g\in G:g\cdot x\in{\cal O}\})=_{df}{\rm{inf}}\{d(h,g):g\cdot x\in 
{\cal O}\}.\]\\

1.17. Lemma. $\varphi^{\cal O}$ is a $G$-mapping. 

Proof. Fix $\bar{g}\in G$. Then for all $h\in G$
\[\varphi^{\cal O}_{\bar{g}\cdot x}(h)=d(h,\{g\in G:g\bar{g}\cdot x\in{\cal O}\},\]
which by definition equals
\[{\rm{inf}}\{d(h,g):g\bar{g}\cdot x\in
{\cal O}\},\]
which by right invariance equals 
\[{\rm{inf}}\{d(h\bar{g},g\bar{g}):g\bar{g}\cdot x\in    
{\cal O}\}={\rm{inf}}\{d(h\bar{g},g'):g'\cdot x\in    
{\cal O}\}\]
\[=_{df}\varphi^{\cal O}_x(h\bar{g})=_{df}(\bar{g}\cdot\varphi^{\cal O}_x)(h),\]
as required.$\Box$\\

Note that for $U$ open, $\varphi^{\cal O}_x(h)>0$ if and only if there is an open set $V\subset G$ 
with $h\in V$ and $\forall g\in V(g\cdot(h\cdot x)\in X\setminus {\cal O})$. 
Thus if $W$ is open in $G$ and ${\cal C}=X\setminus {\cal O}$, then $x\in {\cal C}^{\Delta W}$ if 
and only if there is some $h\in W$ with $\varphi^{\cal O}_x(h)>0$. 
The point here is that this last condition is open in the topology on ${\cal L}(G, d)$.\\

1.18. Notation. For $g$ a point in a metric space $(M,d)$, $q\in{\Bbb Q}^+$ let $B_q^d(g)$ be the ball of 
all points in $M$ with $d(x,g)<q$.\\

\newpage
{\bf \S2. The proof}\\

The rest will be brief.\\

2.1. Lemma. Let $G$ be a Polish group, $(X,\tau)$ a Polish $G$-space, $U_i\subset G$ open, 
$F_i\subset X$ in $\Ubf{\Sigma}^0_2(X,\tau)$, each $i\in{\Bbb N}$. 
Then there is a Polish topology $\tau^*$ on $X$ such that 

(i) $(X,\tau^*)$ is a Polish $G$-space;

(ii) $\tau\subset\tau^*$;

(iii) ${\cal O}\in \tau^*\Rightarrow {\cal O}\in \Ubf{\Sigma}^0_2(X,\tau)$;

(iv) $F_i^{\Delta U_i}\in\tau^*$. 

Proof. By 1.3 it suffices to show that if $({\cal C}_i)_{i\in{\Bbb N}}$ is a 
sequence of closed sets in $(X,\tau)$, then we can find $\tau^*$ satisfying 
(i), (ii), and (iii) above, with ${\cal C}_i^{\Delta W}$ open for each $W\subset G$. 
For notational simplicity let us concentrate on achieving this outcome 
for a single $\tau$-closed set ${\cal C}$; the more general case has an 
exactly similar proof. 

So fix ${\cal C}\subset X$ $\tau$-open, and suppose that ${\cal C}=X\setminus {\cal O}$. 
Let $\hat{X}\subset X\times {\cal L}(G,d)$ be given by 
$\hat{X}=\{(x,\varphi^{\cal O}_x):x\in X\}$. 
By 1.10 the $X\times {\cal L}(G,d)$ is a Polish $G$-space, or rather, 
can be naturally viewed as such, by taking the product of the actions; 
by 
1.7 the map 
\[\pi:X\rightarrow X\times {\cal L}(G,d)\]
\[x\mapsto (x,\varphi^{\cal O}_x)\]
is a $G$-embedding. Let $\tau_0$ be the product topology on $X\times {\cal L}(G,d).$ 

Claim 1. $\hat{X}$ is $\Ubf{\Pi}^0_2(X\times {\cal L}(G,d),\tau_0)$. 

Proof of claim. Let $G_0$ be a countable dense subset of $G$. Then to assert that for 
$x\in X$ and $f\in {\cal L}(G,d)$ we have $f=\varphi_x^{\cal O}$ amounts to 
asserting that 
\[\forall q\in{\Bbb Q}\forall g\in G_0(f(g)\leq q\Leftrightarrow \exists h\in B_q(g) 
(h\cdot x \in {\cal O})).\] 
For any particular $q$ and $g$, the statement that $f(g)\leq q\Leftrightarrow \exists h\in B_q(g)
(h\cdot x \in {\cal O})$ is Boolean combination of $\tau_0$-open sets and hence 
clearly $\Ubf{\Pi}^0_2(X\times {\cal L}(G,d),\tau_0)$, as must therefore be 
the countable intersection over ${\Bbb Q}$ and $G_0$.(Claim$\Box$) 

So $\hat{X}$ is a Polish $G$-space by 1.12, and we have a $G$-embedding 
$\pi:X\rightarrow \hat{X}$, obtained by restricting the map $\pi$ into 
$X\times {\cal L}(G,d)$. Use $\tau_0$ to also denote the topology of 
$X\times {\cal L}(G,d)$ restricted to $\hat{X}$. 

Claim 2. If $W\subset \hat{X}$ is $\tau_0$-open, then $\pi^{-1}(W)$ is 
$\Ubf{\Sigma}^0_2(X,\tau)$. 

Proof of claim. Let $G_0$ be a countable dense subset of $G$, and 
let ${\cal B}=\{B_q(\bar{g}):\bar{g}\in G_0, q\in {\Bbb Q}^+\}$ be a 
countable subbasis of $G$. 
It suffices to prove the claim for a subbasic open set in $\hat{X}$ supported 
on the ${\cal L}(G,d)$ coordinate. 
But for $x\in X$, $g\in G$
\[\varphi_x^{\cal O}(g)>q\Leftrightarrow \exists \bar{q}\in {\Bbb Q}^+
\exists V\in {\cal B}(B_{q+\bar{q}}(g)\subset V\wedge \forall h\in V(h\cdot x\in C)),\]
which is $\Ubf{\Sigma}^0_2(X,\tau)$, while 
\[\varphi_x^{\cal O}(g)<q\Leftrightarrow \exists h\in B_q(g)(h\cdot x\in{\cal O}),\]
which is $\tau$-open.(Claim$\Box$)

Now let $\tau^*$ be the topology consisting of sets of the form $\pi^{-1}(W)$, 
where $W\subset \hat{X}$ is $\tau_0$ open. Since it follows from the 
definition of $\pi$ that $\pi$ is an open mapping, we have $\tau\subset\tau^*$. 
Since $\pi$ is a $G$-embedding, we obtain that $(X,\tau^*)$ is a Polish $G$-space. 
Following the remark after 1.17, ${\cal C}^{\Delta V}$ is open in 
$(X, \tau^*)$ for any open $V\subset G$. 

The more general case, where we consider many ${\cal C}_i$, follows a similar 
proof. Now $\hat{X}$ will be a $\Ubf{\Pi}^0_2$ subset of 
$X\times({\cal L}(G,d))^{\Bbb N}$, where the $i$th copy of ${\cal L}(G,d)$ 
is responsible for making sets of the form ${\cal C}_i^{\Delta V}$ open.$\Box$\\

2.2. Theorem. Let $G$ be a Polish group, $(X,\tau)$ a Polish $G$-space, $U_i\subset G$ open,
$B_i\subset X$ in $\Ubf{\Sigma}^0_{\alpha}(X,\tau)$, for each $i\in{\Bbb N}$ for ${\alpha}$ some 
fixed countable ordinal. Then there is a Polish topology $\tau^*$ on $X$ such that

(i) $(X,\tau^*)$ is a Polish $G$-space;

(ii) $\tau\subset\tau^*$;

(iii) ${\cal O}\in \tau^*\Rightarrow {\cal O}\in \Ubf{\Sigma}^0_{\alpha}(X,\tau)$;

(iv) for each $i$, $B_i^{\Delta U_i}\in\tau^*$.

Proof. By induction on $\alpha$. 

If $\alpha$ is an infinite limit ordinal, 
then we may assume that each $B_i=\bigcup \{B_{i,j}:j\in{\Bbb N}\}$, with 
each $B_{i,j}$ in $\Ubf{\Sigma}^0_{\alpha(j)}(X,\tau)$ some fixed $\alpha(j)<\alpha$. 
Then for ${\cal B}$ a basis for $G$, we can employ our inductive hypothesis to 
find a sequence of 
Polish topologies $(\tau_j)_{j\in{\Bbb N}}$ such that 

(i') each $(X,\tau_j)$ is a Polish $G$-space;

(ii') $\tau\subset\tau_{j}$;

(iii') ${\cal O}\in \tau_{j}\Rightarrow {\cal O}\in 
\Ubf{\Sigma}^0_{\alpha(j)}(X,\tau)$;

(iv') for each $i\in{\Bbb N}$, $W\in{\cal B}$, $B_{i,j}^{\Delta W}\in\tau_j$. 

Let $Y=\Pi\{(X,\tau_j):j\in{\Bbb N}\}$. By 1.13 this is a Polish space. 
Note also that it is in fact a Polish $G$-space in the product action. 
Let $\tau^*$ be the topology generated by taking $\bigcup\{\tau_j:j\in{\Bbb N}\}$ 
a subbasis. Then $(X, \tau^*)$ is homeomorphic as a $G$-space to the 
diagonal $\{\vec x\in Y:\forall i(x(i)=x(i+1))\}$, which is a closed 
invariant subset of $Y$ since each $\tau_i$ includes $\tau$. 
It follows from the definitions that under 
$\pi:X\rightarrow Y$ defined by $(\pi(x))(i)=x$ we have that the 
pullbacks of open sets are all $\Ubf{\Sigma}^0_{\alpha}(X,\tau)$. 
Thus, as before, $\tau^*$ is as required. 

For the inductive step, suppose that $\alpha=\beta+1$ and we have already 
established the theorem for $\beta$. Then we may assume that 
each $B_i=\bigcup\{B_{i,j}:j\in{\Bbb N}\}$, each $B_{i,j}$ is 
$\Ubf{\Pi}^0_{\beta}(X,\tau)$. Let $\tau'$, by the inductive hypothesis, 
be a topology on $X$ so that: 

(i'') $(X,\tau')$ is a Polish $G$-space;

(ii'') $\tau\subset\tau'$;

(iii'') ${\cal O}\in \tau'\Rightarrow {\cal O}\in
\Ubf{\Sigma}^0_{\beta}(X,\tau)$;

(iv'') for each $i\in{\Bbb N}$, $W\in{\cal B}$, $B_{i,j}^{*W}\in\Ubf{\Pi}^0_1(\tau')$.

Thus we can apply 2.1 and find a Polish topology $\tau^*$ such that 

(i$^*$) $(X,\tau^*)$ is a Polish $G$-space;      

(ii$^*$) $\tau\subset\tau'\subset\tau^*$;  

(iii$^*$) ${\cal O}\in \tau^*\Rightarrow {\cal O}\in  
\Ubf{\Sigma}^0_{2}(X,\tau')\subset\Ubf{\Sigma}^0_{\beta+1}(X,\tau)$;    

(iv$^*$) for each $i\in{\Bbb N}$, $W_0,W_1\in{\cal B}$, 
$(B_{i,j}^{*W_0})^{\Delta W_1}\in\Ubf{\Sigma}^0_1(\tau^*)$.

Now we observe that  for $B_i=\bigcup\{B_{i,j}:j\in{\Bbb N}\}$, $W\in{\cal B}$ 
\[B_i^{\Delta W}=\bigcup\{(B_{i,j}^{* W_0})^{\Delta W_1}: W_0, W_1\in{\cal B}, 
W_0\cdot W_1\subset W\},\] 
and thus is open with respect to $\tau^*$.$\Box$\\

2.3. Corollary. Let $G$ be a Polish group, $(X,\tau)$ a Polish $G$-space, 
$B_i\subset X$ a sequence of $G$-invariant sets -- that is, 
$\forall g\in x\in B_i(g\cdot x\in B_i)$. 
Suppose $B_i\in\Ubf{\Sigma}^0_{\alpha}(X,\tau)$. 
Then there is a Polish topology $\tau^*$ on $X$ such that

(i) $(X,\tau^*)$ is a Polish $G$-space;

(ii) $\tau\subset\tau^*$;

(iii) ${\cal O}\in \tau^*\Rightarrow {\cal O}\in \Ubf{\Sigma}^0_{\alpha}(X,\tau)$;

(iv) each $B_i\in\tau^*$.

Proof. By 2.2, since if $U$ is any non-empty open subset of $G$ we have 
that each $B_i^{\Delta U}=B_i$ by invariance.$\Box$

\newpage

{\bf \S3. Some connections}\\

The construction above can be given other tasks. 

For instance, as remarked by Alexander Kechris, 
if $({\cal O}_i)_{i\in{\Bbb N}}$ forms a basis for the topology 
for the Polish $G$-space $(X,\tau)$, then for $d$ a right invariant metric 
we may consider the function 
\[\rho: X\rightarrow ({\cal L}(G,d))^{\Bbb N} \]
\[\rho: x\mapsto (\varphi_x^{{\cal O}_i})_{i\in{\Bbb N}}.\] 
It follows from the proof of 2.1 that $\rho$ is a $G$-map such 
that the pull back of any open set is $\Ubf{\Sigma}^0_2(X, \tau)$ -- 
in other words, $\rho$ is a Baire class 1 function. Since we 
have chosen enough open sets, $\rho$ is in fact a Baire class 1 
$G$-embedding from $(X,\tau)$ to a compact Polish $G$-space. 
Thus we obtain a new proof of a theorem from \cite{beckerkechris} 
that for every Polish group there is a universal Polish $G$-space 
${\cal U}_G$, such that every other Polish $G$ space allows a 
Borel $G$ embedding into ${\cal U}_G$. 

It is known from work of Megrelishvili that in general there may 
exist Polish groups $G$ for which there is a Polish $G$-space that 
allows no continuous $G$-embedding into a compact Polish $G$-space. 
Thus we may choose to view this as a kind of optimal result -- 
one in general cannot hope for a 
compact Polish $G$-space that is universal via continous maps, 
but there does exist a compact Polish $G$-space that is 
universal via Baire class 1 $G$-embeddings. 

Kechris has also shown that in this construction the pointwise 
image of $X$ under $\rho$ is $\Ubf{\Pi}^0_2$ in the product 
topology on $({\cal L}(G,d))^{\Bbb N}$, thereby strengthening 
the sense that this is optimal granted the Megrelishvili counterexample. 

Another application was noted by Ramez Sami, who commented that 
we hereby obtain a generalization of a result from \cite{sami} 
that was initially proved only for the Polish group $S_{\infty}$.\\

3.1. Proposition. Let $G$ be a Polish group, 
$(X,\tau)$ a Polish $G$-space, and suppose that there 
are less than $2^{\aleph_0}$ many orbits. Then any 
invariant $\Ubf{\Pi}^0_{\alpha+1}(X,\tau)$ contains 
a $\Ubf{\Pi}^0_{\alpha+1}(X,\tau)$ orbit. 

Proof. Let $B$ be the invariant set. We can find $\tau^*$ as in 2.2 so 
that $\tau\subset\tau^*$ and $B\in \Ubf{\Pi}^0_{2}(X,\tau^*)$ 
and every $\tau^*$-open set is $\Ubf{\Sigma}^0_{\alpha}(X,\tau)$. 
Since $B$ does not contain $2^{\aleph_0}$ many orbits, 
the orbit equivalence relation is not meager, and hence 
we may find a $\Ubf{\Pi}^0_{2}(X,\tau^*)$ orbit, as in \cite{sami} 
(or see \cite{beckerkechris}). Then this orbit will be 
$\Ubf{\Pi}^0_{\alpha+1}(X,\tau)$, as required.$\Box$\\

Using work from \cite{beckerkechris} one can reduce 
the assumption that there are less than $2^{\aleph_0}$ many orbits 
to the weaker hypothesis that the orbit equivalence relation 
does not {\it Borel reduce} $E_0$, the equivalence relation of eventual 
agreement on countable sequences of integers.\\

{\bf Acknowledgements}  I thank Su Gao and Alexander Kechris, who both 
read through an earlier draft of this note and made  
many helpful suggestions; Kechris in particular 
caught various errors and alerted me to the 
Megrelishvili paper.

6363 MSB 
Mathematics 
UCLA, CA 90095-1555\\

greg@math.ucla.edu

\end{document}